\documentclass[12pt]{article}
\usepackage{amssymb,amsthm,amsmath,latexsym}
\usepackage{url}
\usepackage{color}
\newtheorem{thm}{Theorem}

\newtheorem{lem}[thm]{Lemma}

\theoremstyle{remark}

\newcommand{\FF}{\mathbb{F}}
\newcommand{\ZZ}{\mathbb{Z}}
\newcommand{\RR}{\mathbb{R}}

\begin{document}
\title{
Construction of 
$s$-extremal optimal unimodular lattices in dimension 52
}

\author{
Masaaki Harada\thanks{
Research Center for Pure and Applied Mathematics,
Graduate School of Information Sciences,
Tohoku University, Sendai 980--8579, Japan.
email: {\tt mharada@tohoku.ac.jp}.}
}


\maketitle

\vspace{-1cm}
\begin{center}
{\bf In memory of Masaaki Kitazume}
\end{center}

\begin{abstract}
An $s$-extremal optimal unimodular lattice in dimension
$52$ is constructed for the first time.
This lattice is constructed from a certain self-dual $\FF_5$-code
by Construction A.
In addition, as neighbors of the lattice, 
two more  $s$-extremal optimal unimodular lattices are constructed.
\end{abstract}

\section{Introduction}

A (Euclidean) integral lattice $L$ in dimension $n$
is {\em unimodular} if $L = L^{*}$, where
$L^{*}$ denotes the dual lattice of $L$.
A unimodular lattice with the largest minimum norm among all
unimodular lattices in that dimension is called 
{\em optimal}~\cite{CS98}.
The largest possible minimum norms for dimension $n \le 40$
are given in~\cite[Table 1]{CS98}.
This work was extended in~\cite{NG}.
For example, the largest minimum norms for dimensions $52$ and $54$ are $5$ 
(see~\cite{NG}).

Let $L$ be an odd unimodular lattice in dimension $n$.
The {\em shadow} $S(L)$ of $L$ is defined to be $S(L)= L_0^* \setminus L$,
where $L_0$ denotes the even sublattice of $L$~\cite{CS98}.
It was shown in~\cite{Ga07} that 
$8 \min(L) + 4 \min(S(L)) \le 8+n$, 
unless $n=23$ and $\min(L)=3$, 
where $\min(L)$ and $\min(S(L))$ denote the minimum norms of 
$L$ and $S(L)$, respectively.
A unimodular lattice meeting the above bound with equality is called 
{\em $s$-extremal}.
For minimum norms $\mu \le 4$, $s$-extremal unimodular lattices 
have been widely studied 
(see e.g.\ \cite{Ga07}, \cite{H11} and the references therein).
For minimum norm $\mu \ge 5$, only a few
$s$-extremal unimodular lattices are currently known.
More precisely, 
three nonisomorphic $s$-extremal odd unimodular lattices
in dimension $48$ having minimum norm $5$
were constructed in~\cite{Venkov}, and 
an $s$-extremal optimal unimodular lattice
in dimension $54$ was constructed in~\cite{Ga07}.
The aim of this note is to establish the following theorem.

\begin{thm}\label{main}
There are at least three nonisomorphic 
$s$-extremal optimal unimodular lattices
in dimension $52$.
\end{thm}


For dimension $52$,
no $s$-extremal optimal unimodular lattice was previously known.

%
 

 \section{Preliminaries}\label{Sec:2}

\subsection{Unimodular lattices}
A (Euclidean) lattice $L \subset \RR^n$ 
in dimension $n$
is {\em unimodular} if
$L = L^{*}$, where
the dual lattice $L^{*}$ of $L$ is defined as
$\{ x \in {\RR}^n \mid (x,y) \in \ZZ \text{ for all }
y \in L\}$ under the standard inner product $(x,y)$.
A unimodular lattice $L$ is {\em even} 
if the norm $(x,x)$ of every vector $x$ of $L$ is even,
and {\em odd} otherwise.
An even unimodular lattice in dimension $n$
exists if and only if $n \equiv 0 \pmod 8$, while
an odd  unimodular lattice exists for every dimension.
The {\em minimum norm} $\min(L)$ of 
a unimodular lattice $L$ is the smallest 
norm among all nonzero vectors of $L$.
The {\em kissing number} of $L$ is the number of vectors of minimum norm in $L$.
Two lattices $L$ and $L'$ are {\em isomorphic},
if there is an orthogonal matrix $A$ with
$L' = L \cdot A$, where
$L \cdot A =\{xA \mid x \in L\}$.
The {\em automorphism group} of $L$ is the group of all
orthogonal matrices $A$ with $L = L \cdot A$.

Let $L$ be an odd unimodular lattice in dimension $n$.
The {\em shadow} $S(L)$ of $L$ is defined to be $S(L)= L_0^* \setminus L$,
where $L_0$ denotes the even sublattice of $L$.
Shadows of odd unimodular lattices appeared in~\cite{CS98} 
(see also~\cite[p.~440]{SPLAG}),
and shadows play an important role in the study of odd unimodular lattices.
It was shown in~\cite{Ga07} that 
\begin{equation}\label{eq:B}
8 \min(L) + 4 \min(S(L)) \le 8+n, 
\end{equation}
unless $n=23$ and $\min(L)=3$, 
where $\min(S(L))$ 
denotes the minimum norm of $S(L)$, that is,
the smallest norm among all nonzero vectors of $S(L)$.
An odd unimodular lattice meeting the bound~\eqref{eq:B} with equality is called 
{\em $s$-extremal}.

The {\em theta series} of an odd unimodular lattice $L$ and its shadow $S(L)$
are the formal power series
$\theta_{L}(q) = \sum_{x \in L} q^{(x,x)}$ and
$\theta_{S(L)}(q) = \sum_{x \in S(L)} q^{(x,x)}$, respectively.
Conway and Sloane~\cite{CS98} showed that
when the theta series of an odd unimodular lattice $L$
in dimension $n$
is written as:
\begin{equation}
\label{eq:theta}
 \sum_{j =0}^{\lfloor n/8\rfloor} a_j\theta_3(q)^{n-8j}\Delta_8(q)^j,
\end{equation}
the theta series of the shadow $S(L)$
is written as:
\begin{equation}
\label{eq:theta-S}
\sum_{j=0}^{\lfloor n/8\rfloor}
\frac{(-1)^j}{16^j} a_j\theta_2(q)^{n-8j}\theta_4(q^2)^{8j},
\end{equation}
where
$\Delta_8(q) = q \prod_{m=1}^{\infty} (1 - q^{2m-1})^8(1-q^{4m})^8$,
and $\theta_2(q), \theta_3(q)$, $\theta_4(q)$ are the Jacobi
theta series~\cite{SPLAG}.

\subsection{Self-dual $\FF_5$-codes and Construction A}

Let $\FF_5$ denote the finite field of order $5$.
An {\em $\FF_5$-code} $C$ of length $n$
is a vector subspace of $\FF_5^n$.
The {\em dual} code $C^{\perp}$ of an $\FF_5$-code $C$ of length $n$
is defined as
$
C^{\perp}=
\{x \in \FF_5^n \mid x \cdot y = 0 \text{ for all } y \in C\},
$
where $x \cdot y$ is the standard inner product.
A code $C$ is {\em self-dual} if $C=C^\perp$.

An $n \times n$ negacirculant matrix has the following form:
\[
\left(
\begin{array}{cccccc}
r_0&r_1 &\cdots      &r_{n-2} &r_{n-1} \\
-r_{n-1}&r_0& \cdots &r_{n-3} &r_{n-2} \\
\vdots & \ddots & \ddots&& \vdots\\
-r_2 &        &\ddots & \ddots& r_1\\
-r_1&-r_2& \cdots   &-r_{n-1}&r_0 \\
\end{array}
\right).
\]
Let $A$ and $B$ be $n \times n$ negacirculant matrices.
Consider a matrix of the following form:
\begin{equation} \label{eq:GM}
\left(
\begin{array}{ccc@{}c}
\quad & {\Large I_{2n}} & \quad &
\begin{array}{cc}
A & B \\
-B^T & A^T
\end{array}
\end{array}
\right),
\end{equation}
where $I_{2n}$ denotes the identity matrix of order $2n$,
and $A^T$ denotes the transpose of a matrix $A$.
The codes with generator matrices of the form~\eqref{eq:GM}
are called {\em four-negacirculant}~\cite{4cir}.

Let $C$ be a self-dual $\FF_5$-code of length $n$.
Then the following lattice
\[
A_{5}(C) = \frac{1}{\sqrt{5}}
\{(x_1,\ldots,x_n) \in \ZZ^n \mid
(x_1 \bmod 5,\ldots,x_n \bmod 5)\in C\}
\]
is a unimodular lattice in dimension $n$.
This construction of lattices is well known as {Construction A}.

\section{On $s$-extremal optimal unimodular lattice in 
dimension $52$}

The largest minimum norm among unimodular lattices in dimension $52$
is $5$~\cite{NG}.
By using~\eqref{eq:theta} and \eqref{eq:theta-S},
we determine the possible theta series $\theta^L_{52}(q)$ and
$\theta^S_{52}(q)$ of an optimal odd unimodular lattice 
and its shadow, respectively, as follows.
Since the minimum norm is $5$, 
we have that
\[
a_0=1,
a_1=-104,
a_2=3016,
a_3=-22464,
a_4=19656,
\]
in~\eqref{eq:theta}.
By considering the coefficients of $q$ and $q^3$ in~\eqref{eq:theta-S},
$a_5$ and $a_6$ are written as:
\[
a_5=-256 \alpha,
a_6=1048576 \beta,
\]
by using integers $\alpha$ and $\beta$.
Hence, we have that
\begin{align*}
\theta^L_{52}(q) =&
1 
+ (157248  - 256 \alpha )q^5 
+ (15462720  + 4096 \alpha  + 1048576 \beta )q^6 
\\ & \hspace{2cm}
+ (729181440  - 21504 \alpha  - 41943040 \beta )q^7 
 + \cdots,
\\
\theta^S_{52}(q) =&
\beta q 
+ (\alpha  - 92 \beta) q^3 
+ (314496  - 68 \alpha  + 4134 \beta) q^5 
\\ & \hspace{2cm}
+ (1458362880  + 2226 \alpha  - 120888 \beta) q^7
+ \cdots.
\end{align*}
The first optimal unimodular lattice in dimension $52$ was constructed
in~\cite{Ga04}, and the lattice has theta series 
$\theta^L_{52}(q)$ with $(\alpha,\beta)=(104,0)$.

An easy observation on $\theta^L_{52}(q)$ and $\theta^S_{52}(q)$
gives the following lemma.

\begin{lem}\label{lem:52}
Let $L$ be an optimal odd unimodular lattice in dimension $52$.
Then, the kissing number of $L$ is $157248$
if and only if $L$ is $s$-extremal.
\end{lem}
\begin{proof}
Follows from 
the coefficient of $q^5$ in $\theta^L_{52}(q)$ and
the coefficients of $q$ and $q^3$ in $\theta^S_{52}(q)$.
\end{proof}

For an odd unimodular lattice $L$, 
there are cosets $L_1,L_2,L_3$ of $L_0$ such that
$L_0^* = L_0 \cup L_1 \cup L_2 \cup L_3$, where
$L = L_0  \cup L_2$ and $S(L) = L_1 \cup L_3$.
Two lattices $L$ and $L'$ are {\em neighbors} if
both lattices contain a sublattice of index $2$
in common.
If $L$ is an odd unimodular lattice in dimension divisible by
$4$, then there are two unimodular lattices
containing $L_0$,
which are rather than $L$,
namely, $L_0 \cup L_1$ and $L_0 \cup L_3$ (see~\cite{DHS}).
We denote the two unimodular neighbors by
\[
N_1(L)=L_0 \cup L_1 \text{ and } N_2(L)=L_0 \cup L_3.
\]

\begin{lem}\label{lem}
If $L$ is an $s$-extremal optimal unimodular lattice in dimension $52$, then
the lattices $N_1(L)$ and $N_2(L)$ are also
$s$-extremal optimal unimodular lattices in that dimension.
\end{lem}
\begin{proof}
Since $L$ has shadow of minimum norm $5$,
both $N_1(L)$ and $N_2(L)$ have
minimum norm $5$.
The shadow of $N_1(L)$ is
$L_2 \cup L_3$ and
the shadow of $N_2(L)$ is
$L_2 \cup L_1$.
Hence, 
both $N_1(L)$ and $N_2(L)$ 
have shadows of minimum norm $5$, that is, $s$-extremal.
\end{proof}


\section{Three examples of $s$-extremal optimal unimodular lattices in 
dimension $52$}

By considering four-negacirculant codes, we found a certain
self-dual $\FF_5$-code $C_{52}$ of length $52$.
The self-duality of $C_{52}$ was verified by the {\sc Magma} function
{\tt IsSelfDual} (see~\cite{Magma} for {\sc Magma}).
The code $C_{52}$ has generator matrix~\eqref{eq:GM},
where the first rows of the negacirculant matrices
$A$ and $B$ in~\eqref{eq:GM} are as follows:
\[
(1,2,3,3,4,3,2,0,3,1,0,3,1) \text{ and }
(4,4,0,3,1,0,4,2,0,1,3,3,1),
\]
respectively.
In addition, we verified that
$A_5(C_{52})$ has minimum norm $5$ and kissing number $157248$.
The minimum norm and the kissing number of 
$A_5(C_{52})$ were calculated
by the {\sc Magma} functions {\tt Minimum}
and {\tt KissingNumber}, respectively.
By Lemma~\ref{lem:52}, 
$A_5(C_{52})$ is an $s$-extremal optimal
unimodular lattice in dimension $52$.
By Lemma~\ref{lem}, 
$N_1(A_5(C_{52}))$ and $N_2(A_5(C_{52}))$
are $s$-extremal optimal
unimodular lattices in dimension $52$.
In addition,
we verified by the {\sc Magma} function {\tt IsIsomorphic} that 
$A_5(C_{52})$,  $N_1(A_5(C_{52}))$ and $N_2(A_5(C_{52}))$
are nonisomorphic.
Hence, we have established  Theorem~\ref{main}.
Finally, we verified that 
$A_5(C_{52})$,  $N_1(A_5(C_{52}))$ and $N_2(A_5(C_{52}))$
have automorphism groups of order $52$.
This was verified by the {\sc Magma} function {\tt AutomorphismGroup}.

\bigskip
\noindent
{\bf Acknowledgment.}
This work was supported by JSPS KAKENHI Grant Number 19H01802.



\end{document}